\documentclass{article}
\usepackage{amsmath,amssymb}

\newtheorem{thm}{Theorem}[section]

\newtheorem{prop}[thm]{Proposition}

\let\a=\alpha    
    
  \let\n=\nu

\let\C=\Chi

\def\nn{\nonumber} \def\bd{\begin{document}} \def\ed{\end{document}}
\def\ds{\documentstyle} \let\fr=\frac \let\bl=\bigl \let\br=\bigr
\let\Br=\Bigr \let\Bl=\Bigl 
\let\bm=\bibitem
\let\na=\nabla
\let\pa=\partial \let\ov=\overline 
\def\ba{\begin{array}}
\def\ea{\end{array}}
\def\ft#1#2{{\textstyle{\frac{\scriptstyle #1}{\scriptstyle #2}}}}
\def\fft#1#2{\frac{#1}{#2}}
\def\del{\partial}
\def\vp{\varphi}
\def\st#1{{\scriptstyle #1}}
\def\sst#1{{\scriptscriptstyle #1}}

\def\oneone{\rlap 1\mkern4mu{\rm l}}
\def\td{\tilde}
\def\wtd{\widetilde}
\def\ie{\rm i.e.\ }
\def\dalemb#1#2{{\vbox{\hrule height .#2pt
        \hbox{\vrule width.#2pt height#1pt \kern#1pt
                \vrule width.#2pt}
        \hrule height.#2pt}}}
\def\square{\mathord{\dalemb{6.8}{7}\hbox{\hskip1pt}}}

\def\cramp{\medmuskip = 2mu plus 1mu minus 2mu}
\def\cramper{\medmuskip = 2mu plus 1mu minus 2mu}
\def\crampest{\medmuskip = 1mu plus 1mu minus 1mu}
\def\uncramp{\medmuskip = 4mu plus 2mu minus 4mu}

\newcommand{\ho}[1]{$\, ^{#1}$}
\newcommand{\hoch}[1]{$\, ^{#1}$}
\newcommand{\ra}{\rightarrow}
\newcommand{\lra}{\longrightarrow}
\newcommand{\Lra}{\Leftrightarrow}
\newcommand{\ap}{\alpha^\prime}
\newcommand{\bp}{\tilde \beta^\prime}
\newcommand{\tr}{{\rm tr} }
\newcommand{\Tr}{{\rm Tr} } 
\def\0{{\sst{(0)}}}
\def\1{{\sst{(1)}}}
\def\2{{\sst{(2)}}}
\def\3{{\sst{(3)}}}
\def\4{{\sst{(4)}}}
\def\5{{\sst{(5)}}}
\def\6{{\sst{(6)}}}
\def\7{{\sst{(7)}}}
\def\8{{\sst{(8)}}}
\def\n{{\sst{(n)}}}
\def\cA{{{\cal A}}}
\def\cF{{{\cal F}}}
\def\tV{\widetilde V}
\def\tW{\widetilde W}
\def\tH{\widetilde H}
\def\tE{\widetilde E}
\def\tF{\widetilde F}
\def\tA{\widetilde A}
\def\im{{{\rm i}}}
\def\jm{{{\rm j}}}
\def\km{{{\rm k}}}

\def\tY{{{\wtd Y}}}
\def\ep{{\epsilon}}
\def\vep{{\varepsilon}}
\def\R{\rlap{\rm I}\mkern3mu{\rm R}}
\def\bD{{{\bar D}}}
\def\R{{{\mathbb R}}}
\def\C{{{\mathbb C}}}
\def\H{{{\mathbb H}}}
\def\CP{{{\mathbb C}{\mathbb P}}}
\def\RP{{{\mathbb R}{\mathbb P}}}
\def\Z{{{\mathbb Z}}}
\def\bA{{{\mathbb A}}}
\def\bB{{{\mathbb B}}}

\newcommand{\NP}{Nucl. Phys. }
\newcommand{\tamphys}{\it Center for Theoretical Physics\\
Texas A\&M University, College Station, TX 77843, USA}
\newcommand{\umich}{\it Michigan Center for Theoretical Physics\\
University of Michigan, Ann Arbor, Michigan 48109, USA}
\newcommand{\upenn}{\it Department of Physics and Astronomy\\
University of Pennsylvania, Philadelphia,  PA 19104, USA}
\newcommand{\SISSA}{\it  SISSA-ISAS and INFN, Sezione di Trieste\\
Via Beirut 2-4, I-34013, Trieste, Italy}

\newcommand{\ihp}{\it Institut Henri Poincar\'e\\
  11 rue Pierre et Marie Curie, F 75231 Paris Cedex 05}

\newcommand{\damtp}{\it DAMTP, Centre for Mathematical Sciences,
 Cambridge University, Wilberforce Road, Cambridge CB3 OWA, UK}

\newcommand{\auth}{M. Cveti\v{c}\hoch{\dagger}, G.W. Gibbons\hoch{\sharp}, 
H. L\"u\hoch{\star} and C.N. Pope\hoch{\ddagger}}

\begin{document}
\begin{flushright}
\hfill{DAMTP-2001-39}\ \ \ {CTP TAMU-15/01}\ \ \ {UPR-937-T}\ \ \
{MCTP-01-23}\\ 
{May 2001}\\
{math.DG/0105119}
\end{flushright}

\vspace{15pt}

\begin{center}
{ \large {\bf New Cohomogeneity One Metrics With Spin(7) Holonomy}}

\vspace{15pt}
\auth

\vspace{5pt}
{\hoch{\dagger}\upenn}

\vspace{5pt}
{\hoch{\sharp}\damtp}

%%\vspace{5pt}
%%{\hoch{\dagger} \it Department of Physics and Astronomy, Rutgers University,
%%Piscataway, NJ 08855}

\vspace{4pt}
{\hoch{\star}\umich}

\vspace{4pt}
{\hoch{\ddagger}\tamphys}

%\vspace{5pt}
%{\hoch{\dagger,\sharp,\ddagger}\ihp}

\vspace{10pt}

\underline{ABSTRACT}
\end{center}

   We construct new explicit non-singular metrics that are complete on
non-compact Riemannian 8-manifolds with holonomy Spin(7).  One such
metric, which we denote by $\bA_8$, is complete and non-singular on
$\R^8$.  The other complete metrics are defined on manifolds with the
topology of the bundle of chiral spinors over $S^4$, and are denoted
by $\bB_8^+$, $\bB_8^-$ and $\bB_8$.  The metrics on $\bB_8^+$ and
$\bB_8^-$ occur in families with a non-trivial parameter.  The metric
on $\bB_8$ arises for a limiting value of this parameter, and locally
this metric is the same as the one for $\bA_8$.  The new Spin(7)
metrics are asymptotically locally conical (ALC): near infinity they
approach a circle bundle with fibres of constant length over a cone
whose base is the squashed Einstein metric on $\CP^3$.  We construct
the covariantly-constant spinor and calibrating 4-form.  We also
obtain an $L^2$-normalisable harmonic 4-form for the $\bA_8$ manifold,
and two such 4-forms (of opposite dualities) for the $\bB_8$ manifold.

%{\vfill\leftline{}\vfill
%\vskip 5pt
%\footnoterule
%{\footnotesize \hoch{1} Research supported in part by DOE grant
%DE-FG02-95ER40893 and NATO grant 976951. \vskip -12pt} \vskip 14pt
%{\footnotesize \hoch{2} Research supported in full by DOE grant
%DE-FG02-95ER40899 \vskip -12pt} \vskip 14pt
%{\footnotesize  \hoch{3} Research supported in part by DOE
%grant DE-FG03-95ER40917.\vskip  -12pt}}

%\baselineskip=24pt
\pagebreak
\setcounter{page}{1}

%\tableofcontents
\vfill\eject

\section{Introduction}

    Few explicit examples of complete non-compact manifolds admitting
Ricci-flat metrics with the exceptional holonomies $G_2$ in seven
dimensions or Spin(7) in eight dimensions are known.  Three
asymptotically conical examples have been found in $D=7$, for
manifolds with the topology of the bundle of self-dual 2-forms on
$S^4$ or $\CP^2$, and the spin bundle of $S^3$
\cite{brysal,gibpagpop}.  In $D=8$ the only Spin(7) example that was
known was defined on the chiral spin bundle of $S^4$
\cite{brysal,gibpagpop}.

    In this paper we give a construction of new eight-dimensional
metrics of Spin(7) holonomy, and show that among these are examples
that are complete on two different non-compact manifolds.  It
represents and elaboration and simplification of the original
construction in \cite{cglps7}.  The new metrics are all asymptotically
locally conical (ALC), locally approaching $\R\times S^1\times \CP^3$.
The radius of the $S^1$ is asymptotically constant, so the metric
approaches an $S^1$ bundle over a cone with base $\CP^3$.  However,
the Einstein metric on the $\CP^3$ at the base of the cone is not the
Fubini-Study metric, but instead the ``squashed'' metric described as
an $S^2$ bundle over $S^4$.  The new solutions can have very different
short-distance behaviours, with one approaching flat $\R^8$ whilst all
the others approach $\R^4\times S^4$ locally.  The global topology is
that of $\R^8$ in the first case and the bundle of positive (or
negative) chirality spinors over $S^4$ for the others.  An intriguing
feature of two of the new metrics, one on each of the inequivalent
topologies, is that in local coordinates the metrics are identical.
Globally, the metric is complete on a manifold of $\R^8$ topology if
the radial coordinate $r$ is taken to be positive, whilst in the
region with negative $r$ it is instead complete on the manifold
${\mathbb S}(S^4)$ of the bundle of chiral spinors over $S^4$.  We
shall denote the new Spin(7) manifold with $\R^8$ topology by $\bA_8$,
and the new related manifold with ${\mathbb S}(S^4)$ topology by
$\bB_8$.  The more general classes of new manifolds with the topology
of the chiral spin bundle over $S^4$ will be denoted by $\bB_8^+$ and
$\bB_8^-$.

    Our construction is a generalisation of the one that leads to the
previously-known metric of Spin(7) holonomy.  That example is given, 
in local coordinates, by \cite{brysal,gibpagpop}
%%%%%
\begin{equation}
ds_8^2 = \Big[1- \Big(\fft{r_0}{r}\Big)^{10/3}\Big]^{-1} \, dr^2 
         + \ft{9}{100}\, r^2\, \Big[1- 
               \Big(\fft{r_0}{r}\Big)^{10/3}\Big]\,
         h_i^2 + \ft{9}{20} r^2\, d\Omega_4^2\,,\label{spin7metric}
\end{equation}
%%%%%
where
%%%%%
\begin{equation}
h_i\equiv \sigma_i - A^i\,,
\end{equation}
%%%%%
the $\sigma_i$ are left-invariant 1-forms on $SU(2)$, $d\Omega_4^2$ is
the metric on the unit 4-sphere, and $A^i$ is the potential of the
BPST $SU(2)$
Yang-Mills instanton on $S^4$.  The $\sigma_i$ can be written in terms
of Euler angles as
%%%%%
\crampest
\begin{equation}
\sigma_1 = \cos\psi\, d\theta + \sin\psi\, \sin\theta\,
d\varphi\,,\quad
\sigma_2 = -\sin\psi\, d\theta + \cos\psi\, \sin\theta\,
d\varphi\,,\quad
\sigma_3 = d\psi + \cos\theta\, d\varphi\,.\label{1forms}
\end{equation}
\uncramp
%%%%%
The principal orbits are $S^7$, viewed as an $S^3$ bundle over $S^4$.
The solution (\ref{spin7metric}) is asymptotic to a cone over the
``squashed'' Einstein 7-sphere, and it approaches $\R^4\times S^4$
locally at short distance (\ie $r\approx\ell$).  Globally the manifold
has the same topology ${\mathbb S}(S^4)$, the bundle of chiral spinors
over $S^4$, as the new Spin(7) manifolds $\bB_8$ and $\bB_8^\pm$ 
that we obtain in this paper.

\section{Einstein equation and first integrals for Spin(7) metrics}

   The generalisation that we shall consider involves allowing the
$S^3$ fibres of the previous construction themselves to be
``squashed.''  In particular, this encompasses the possibility of
having an asymptotic structure of the ``Taub-NUT type,'' in which the
$U(1)$ fibres in a description of $S^3$ as a $U(1)$ bundle over $S^2$
approach constant length while the radius of the $S^2$ grows linearly.

   A convenient way to parameterise the metric is by first introducing the
left-invariant 1-forms $L_{AB}$ for the
group manifold $SO(5)$.  These satisfy $L_{AB}=-L_{BA}$, and
%%%%%
\begin{equation}
dL_{AB} = L_{AC}\wedge L_{CB}\,.
\end{equation}
%%%%%
The 7-sphere is then given by the coset $SO(5)/SU(2)_L$, where we take
the obvious $SO(4)$ subgroup of $SO(5)$, and write it (locally) as
$SU(2)_L\times SU(2)_R$.
    If we take the indices $A$ and $B$ in $L_{AB}$ to range over the
values $0\le A\le 4$, and split them as $A=(a,4)$, with $0\le a\le 3$,
then the $SO(4)$ subgroup is given by $L_{ab}$.  This is decomposed as
$SU(2)_L\times SU(2)_R$, with the two sets of $SU(2)$ 1-forms given by
the self-dual and anti-self-dual combinations:
%%%%%
\begin{equation}
R_i = \ft12(L_{0i} + \ft12\ep_{ijk}\, L_{jk})\,,\qquad
L_i = \ft12(L_{0i} - \ft12\ep_{ijk}\, L_{jk})\,,
\end{equation}
%%%%%
where $1\le i\le 3$.  Thus the seven 1-forms in the $S^7$ coset will
be 
%%%%%
\begin{equation}
R_1\,,\quad R_2\,,\quad R_3\,,\quad P_a \equiv L_{a4}\,.
\end{equation}
%%%%%
It is straightforward to establish that
%%%%%
\begin{eqnarray}
dP_0 &=& (R_1+ L_1)\wedge P_1 + (R_2+L_2)\wedge P_2 + (R_3+L_3)\wedge
      P_3\,,\nn\\
dP_1 &=& -(R_1+ L_1)\wedge P_0 - (R_2-L_2)\wedge P_3 + (R_3-L_3)\wedge
      P_2\,,\nn\\
dP_2 &=& (R_1- L_1)\wedge P_3 - (R_2+L_2)\wedge P_0 - (R_3-L_3)\wedge
      P_1\,,\nn\\
dP_3 &=& -(R_1- L_1)\wedge P_2 + (R_2-L_2)\wedge P_1 - (R_3+L_3)\wedge
      P_0\,,\nn\\
dR_1 &=& -2 R_2\wedge R_3 - \ft12 (P_0\wedge P_1 + P_2\wedge
P_3)\,,\nn\\
dR_2 &=& -2 R_3\wedge R_1 - \ft12 (P_0\wedge P_2 + P_3\wedge P_1)\,,\nn\\
dR_3 &=& -2 R_1\wedge R_2 - \ft12 (P_0\wedge P_3 + P_1\wedge P_2)\,.
\label{so5d}
\end{eqnarray}
%%%%%
In terms of these left-invariant 1-forms, we can write the ansatz for
the more general metrics of Spin(7) holonomy on
the $\R^4$ bundle over $S^4$ as
%%%%%
\begin{equation}
ds_8^2 = dt^2 + 4a^2\, (R_1^2+R_2^2) + 4b^2\, R_3^2  + c^2\,
P_a^2\,.\label{8metans1}
\end{equation}
%%%%%
We shall work with an orthonormal frame bundle $e^A$ defined by
%%%%%
\begin{equation}
e^8=dt\,,\quad e^{\hat 1} = 2 a\, R_1\,,\quad e^{\hat 2} = 2 a\, R_2\,,\quad
e^{\hat 3} = 2 b\, R_3\,,\quad e^a = c\, R_a\,,
\end{equation}
%%%%% 
where we take the index $a$ to range over $0\le a \le 3$.

   The factors of 4 in the terms involving $a^2$ and $b^2$ in
(\ref{8metans1}) are included for consistency with the conventions in
\cite{cglps7}.  In that paper, the metric was written as
%%%%%
\begin{equation}
d\hat s_8^2 = dt^2 + a^2\, (D\mu^i)^2 + b^2\, \sigma^2 + c^2\,
d\Omega_4^2\,,\label{8metans2}
\end{equation}
%%%%%
where
\begin{equation}
D\mu^i\equiv d\mu^i +\ep_{ijk}\, A^j\, \mu^k\,,\quad
\sigma \equiv d\varphi + \cA\,,\quad \cA \equiv \cos\theta\, d\psi -
\mu^i\, A^i\,,\label{kkvector}
\end{equation}
%%%%%
and $\mu_i$ are coordinates on $\R^3$ subject to the constraint
$\mu_i\, \mu_i=1$ that defines the unit 2-sphere.  In terms of the
left-invariant 1-forms $R_i$ and $P_a$ of this paper, we have the
correspondences
%%%%%
\begin{equation}
 (D\mu^i)^2 = 4(R_1^2+R_2^2)\,,\quad \sigma= 2R_3\,,\quad d\Omega_4^2=
P_a^2\,,\label{corr}
\end{equation}
%%%%%
together with $dA^i + \ft12 \ep_{ijk}\, A^j\wedge A^k = J^i$.  

   We can derive first-order differential equations for the three
metric functions $a$, $b$ and $c$, which will imply Spin(7) holonomy, by
requiring the existence of a closed self-dual 4-form that has the
symmetries of the octonionic structure constants.
It is straightforward to see that an appropriate 4-form, invariant
under the isometries of the metric and with
the required octonionic structure, is given by
%%%%%
\begin{equation}
\Phi_\4 = -e^0\wedge e^1\wedge e^2\wedge e^3
 -e^8\wedge e^{\hat 1} \wedge e^{\hat 2}\wedge e^{\hat 3}
+ \ft12 \ep_{ijk}\, e^{\hat i}\wedge e^{\hat j}\wedge
\hat J^k +
e^8\wedge e^{\hat i}\wedge \hat J^i\,,\label{calib}
\end{equation}
%%%%%
with $\hat J^i\equiv c^2\, J^i= \ft12 J^i_{ab}\, e^a\wedge e^b$, $1\le
i\le 3$, where
%%%%%
\begin{equation}
J^i\equiv P_0\wedge P_i + \ft12 \ep_{ijk}\, P_j\wedge P_k
\end{equation}
%%%%%
are the three self-dual quaternionic-K\"ahler 2-forms on the unit $S^4$
metric $P_a^2$.

   From (\ref{so5d}) we see that $d\Phi_\4=0$ implies the
first-order equations
%%%%%
\begin{equation}
\dot a =1 - \fft{b}{2a} - \fft{a^2}{c^2}\,,\qquad
\dot b  = \fft{b^2}{2a^2} - \fft{b^2}{c^2}\,,\qquad \dot c = \fft{a}{c}
+ \fft{b}{2c}\,.\label{firstorder}
\end{equation}
%%%%%
These equations imply that the metric (\ref{8metans1}) has Spin(7)
holonomy, and also, therefore, that it is Ricci-flat.

   It is sometimes convenient to express the equations for $a$, $b$
and $c$ as a Lagrangian system.  We find that the equations for
Ricci-flatness of (\ref{8metans1}) can be derived by requiring that
$L\equiv T-V$ be stationary with respect to variations of $\a$,
$\beta$ and $\gamma$, where
%%%%%
\begin{eqnarray}
T &=& 2{\a'}^2 + 12 {\gamma'}^2 + 4\a'\, \beta' + 8 \beta'\, \gamma' +
    16 \a'\, \gamma'\,,\nn\\
V &=& \ft12  b^2\, c^4\, (4a^6 + 2a^4\, b^2 - 24 a^4 c^2 -4 a^2 c^4 +
    b^2 \, c^4)\,,\label{tveq}
\end{eqnarray}
%%%%%
together with the constraint $T+V=0$.  Here a prime denotes a derivative
with respect to a new radial variable $\rho$, defined by
$dt=a^2\, b\, c^4\, d\rho$, and we have also defined $\a=\log a$,
$\beta = \log b$,  $\gamma=\log c$.

  We find that the potential $V$ can be derived from a superpotential,
which we denote by $W$.
Writing $T=\ft12 g_{ij}\, (d\a^i/d\rho)\, (d\a^j/d\rho)$,
where $\a^i=(\a,\beta,\gamma)$,  we have
$V=-\ft12 g^{ij}\, (\del W/\del\a^i)\,
(\del W/\del\a^j)$, where
%%%%%
\begin{equation}
W = b\, c^2\, (4a^3 + 2 a^2\, b + 4 a\, c^2-b\, c^2)\,.
\end{equation}
%%%%%
From this we can obtain the first-order equations $d\a^i/d\eta = g^{ij}\,
\del W/\del\a^j$.  Expressed back in terms of the original radial variable
$t$ introduced in (\ref{8metans1}), these equations are precisely those
given in (\ref{firstorder}).

   To summarise, we have the following

\begin{prop} Metrics of the form (\ref{8metans1}) are Ricci-flat if and
only if $a$, $b$ and $c$ satisfy the Euler-Lagrange equations
following from (\ref{tveq}).  The gradient flow defined by
(\ref{firstorder}) provides three first integrals which solve these 
second-order Euler-Lagrange equations.\label{ricciprop}
\end{prop}

    Before proceeding to find new solutions to these first-order
equations, we can first verify that the
previous Spin(7) metric (\ref{spin7metric}) is indeed a solution.
Also, we may observe that one of the seven-dimensional metrics of
$G_2$ holonomy has principal orbits that are $\CP^3$, viewed as an
$S^2$ bundle over $S^4$, and is given, in local coordinates,  
by \cite{brysal,gibpagpop}
%%%%%
\begin{equation}
ds_7^2 = (1-\fft{\ell^4}{r^4})^{-1}\, dr^2 + r^2\, 
(1-\fft{\ell^4}{r^4}) \,(R_1^2+R_2^2) + \ft12 r^2\, P_a^2\,.
\end{equation}
%%%%%
This therefore gives a solution in $D=8$ of
the form $d\hat s_8^2 = ds_7^2 + d\varphi^2$, and it can be described
within the framework of our first-order equations (\ref{firstorder})
by first rescaling $b\longrightarrow\lambda\, b$, and then sending
$\lambda$ to zero, so that $b=$constant is allowed as a solution.

   One can also see the specialisations to the previous results
described above at the level of the first-order equations themselves.
Setting $a=b$ gives a consistent truncation of (\ref{firstorder}),
yielding $\dot a= \ft12 a^2\, c^{-2}$, $\dot c= \ft32 a\, c^{-1}$,
which are indeed the first-order equations for the original Spin(7)
metrics.  On the other hand, sending $b\longrightarrow 0$ in
(\ref{firstorder}) yields a consistent truncation to $\dot a =1-a^2\,
c^{-2}$, $\dot c= a\, c^{-1}$, which are the first-order equations for
the metrics of $G_2$ holonomy whose principal orbits are $S^2$ bundles
over $S^4$.  (The first-order equations for these two cases can be
found, for example, in \cite{cglpcal}.)  Also, we may note that a
special solution arises if we set $b=-a$, which then implies $a=-b=\pm
c= \ft12 t$.  This is flat space.

   Another specialisation of the metric ansatz (\ref{8metans1}) that makes
contact with previous results is to set $a=c$, in which case the $S^2$
bundle over $S^4$ becomes precisely the usual $\CP^3$ Einstein
manifold, with its $SU(4)$-invariant metric.  This is incompatible
with the first-order equations (\ref{firstorder}), but it is easily
verified that it is consistent with the second-order Einstein
equations following from (\ref{tveq}).  Solutions to these
second-order equations then include the 8-dimensional Taub-NUT and
Taub-Bolt metrics.  The incompatibility with the first-order equations
is understandable, since the Taub-NUT and Taub-Bolt 8-metrics do not
have special holonomy.  Another previously-seen solution of the
second-order equations with $a=c$ is the Ricci-flat K\"ahler metric on
the complex line-bundle over $\CP^3$.  Although this can arise from a
first-order system, it is an inequivalent one that is not related to a
specialisation of (\ref{firstorder}).  Its superpotential is $W=2 a^6
+ 6 a^4\, b^2$ \cite{cglp1}, with $T$, $V$ and $g_{ij}$ following from
setting $a=c$ in (\ref{tveq}).  (Other examples of this kind of
phenomenon were exhibited recently in \cite{cglpcal}.)

\section{General solution of the gradient flow}

\subsection{The general solution; local analysis}

    In order to obtain new solutions of the first-order equations 
(\ref{firstorder}) we first introduce a new radial 
coordinate $r$, defined in terms of $t$ by $dr = b\, dt$.  After also defining
$f\equiv c^2$, we find by taking further derivatives of the first-order
equations (\ref{firstorder}) that $f$ must satisfy the third-order 
equation
%%%%%
\begin{equation}
2 f^2\, f''' + 2 f\,(f'-3)\, f'' - (f'+1)(f'-1)(f'-3)=0\,,
\label{feq}
\end{equation}
%%%%%
which can be expressed in the ``factorised'' form 
%%%%%
\begin{equation}
f\, Q' - (f'+1)\, Q=0\,,\label{factored}
\end{equation}
%%%%%
where $Q\equiv 2f\, W' + (f'-3)\, W$ and $W\equiv f'-1$. 
The remaining metric functions are then given by solving
%%%%%
\begin{equation}
a'=\fft{f'-2}{2a} -\fft{(f'-1)\, a}{2f}\,,\qquad
b = \fft{2 a}{(f'-1)}\,.\label{absol}
\end{equation}
%%%%%
Naively there now appear to be four constants of integration in total rather
than the expected three, but the extra one is eliminated by 
substituting the solutions back into (\ref{firstorder}).  In fact for
a generic solution, where $Q$ itself is non-zero, the solution for
$a$, and hence for $b$, can be written entirely algebraically in terms of $f$,
with
%%%%%
\begin{equation}
a^2= \fft{(f'-1)(f'-3)\, f}{Q}\,,\qquad b = \fft{2a}{(f'-1)}\,.
\label{absol2}
\end{equation}
%%%%%
Thus for a solution where $Q\ne0$ the three integration constants for
the first-order system (\ref{firstorder}) are simply the three
integration constants for the third-order equation (\ref{feq}), and no
further substitution back into (\ref{feq}) is necessary.  As we shall
see below, $Q$ is non-vanishing for all but one degenerate solution of
(\ref{feq}).  Note that two of the three constants of integration are 
``trivial,'' corresponding to a constant shift and rescaling of the
radial coordinate.

   The general solution to equation (\ref{feq}) may be obtained as
follows.  First, introduce a new radial variable $\rho$, and a
function $\gamma(\rho)$, defined by
%%%%%
\begin{equation}
f(r) = \exp\Big[ -\int^\rho \fft{\td\rho\, d\td\rho}{\gamma(\td\rho)}\Big]
\,,\qquad \fft{df}{dr} = \rho\,,\label{subs}
\end{equation}
%%%%%
implying also that
%%%%%
\begin{equation}
\fft{d^2f}{dr^2} = -\fft{\gamma}{f}\,, \qquad
\fft{d^3 f}{dr^3} = \fft{1}{f^2}\, \Big(\fft{\rho}{\gamma} + \gamma\,
\fft{d\gamma}{d\rho}\Big)\,.
\end{equation}
%%%%%
(We assume here, and
in the rest of this subsection, that $df/dr$ is not a constant, and so
$\rho$ is a good radial variable.  The special cases where $df/dr$ is
a constant are included in the discussion in section 3.2.)
Equation (\ref{feq}) now reduces to
%%%%%
\begin{equation}
2\gamma\, \fft{d\gamma}{d\rho} + 6\gamma = (1-\rho^2)(3-\rho)\,.
\label{gammaeq}
\end{equation}
%%%%%
The further replacement of $\gamma$ by $z$, and $\rho$ by $v$, defined by
%%%%%
\begin{equation}
z\equiv \fft{(1-\rho)^2}{2(1-\rho-\gamma)}\,,\qquad 
v\equiv \rho-1\label{zvdef}
\end{equation}
%%%%%
turns (\ref{gammaeq}) into
%%%%%
\begin{equation}
2z\, (1-z^2)\, \fft{dv}{dz} = v+ 2z\,.\label{vde}
\end{equation}
%%%%% 
The solution to this equation can be written in terms of the
hypergeometric function as 
%%%%%
\begin{equation}
v = \fft{2k\, \sqrt{z}}{(1-z^2)^{1/4}} -2z\,\,
_2F_1[1,\ft12;\ft54; 1-z^2]\,.\label{vres}
\end{equation}
%%%%% 

    Retracing the steps of the various redefinitions, we see that by
using $z$ as the radial variable the general solution (with $df/dr$
not a constant) for the Ricci-flat metric can be written as
%%%%%
\begin{equation}
a^2 = \fft{(v-2)\, z\, f}{(1+z)\, v}\,,\qquad b = \fft{2
a}{v}\,,\qquad c^2\equiv f = \Big(\fft{1+z}{1-z}\Big)^{1/2}\,
\exp\Big[\int^z \fft{d\td z}{v(\td z)\, (1-{\td z}^2)}\Big]\,.\label{general} 
\end{equation}
%%%%%
The coordinate $r$ is given in terms of $z$ by
%%%%%
\begin{equation}
dr=\fft{f\, dz}{v(z)\, (1-z^2)}\,,
\end{equation}
%%%%%
and so we have the following

\begin{prop} The local Ricci-flat metrics arising from the gradient
flow (\ref{firstorder}) are given by 
%%%%%
\begin{equation}
ds_8^2 = \fft{v\, f\, dz^2}{4z\, (1-z^2)(1-z)\, (v-2)} +
\fft{4(v-2)\,z\, f}{(1+z)\, v}\, (R_1^2+R_2^2) +
  \fft{16(v-2)\, z\, f}{(1+z)\, v^3}\, R_3^2 + f\, P_a^2\,,
\label{general2}
\end{equation}
%%%%% 
with $v$ defined by (\ref{vres}) and $f$ defined by (\ref{general}).
\label{gensolprop}
\end{prop}

    Note also that $k$ in the solution (\ref{vres}) for $v(z)$ is the
non-trivial third constant of integration of the original first-order
system (\ref{firstorder}).  In obtaining this
general solution we have assumed that $Q$ is non-zero, so that $a$ can
be obtained using (\ref{absol2}).  In fact if $Q$ is zero it can
easily be seen that unless in addition $df/dr$ is a constant (which we
have excluded from the analysis in this subsection), then after using
(\ref{absol}) and substituting back into (\ref{firstorder}), the
metric functions $a$, $b$ and $c$ would all vanish.  Thus the only
additional solutions to (\ref{firstorder}), other than those described
by (\ref{general2}), are those with $df/dr=$constant, and these are
included in the discussion in section 3.2.

\subsection{Special globally-defined solutions}

   As we shall show later, the general Ricci-flat metrics obtained in
section 3.1 include one-parameter families of examples
that are complete on manifolds with the topology of the bundle of
chiral spinors over $S^4$.  The parameter in question is a non-trivial
one, as opposed to the two trivial parameters associated with a
constant shift and scaling of the radial coordinate.  Before
discussing these families of complete metrics, we shall first discuss
some simple solutions of the first-order equations (\ref{firstorder}).
It is easier to discuss these in terms of the original radial variable
$r$ used in writing (\ref{feq}).  After absorbing a trivial constant
shift of the radial variable we can write down three elementary
solutions of (\ref{feq}), namely
%%%%%
\begin{equation}
f=-r\,,\qquad f= 3r\,,\qquad f= r + \fft{r^2}{2\ell^2}\,,
\label{elsol}
\end{equation}
%%%%%
where $\ell$ is a constant. The first two solutions here are of the
type where $df/dr$ is a constant, which were excluded in the general
analysis section 3.1.
\medskip

\noindent
$\bullet$   The solution with $f=-r$ has $Q=8$ and so we can use
(\ref{absol2}), to find
%%%%%
\begin{equation}
a^2=-r\qquad b=-a\,,\qquad c^2=-r\,.
\end{equation}
%%%%%
It follows from (\ref{8metans1}) that the metric in this case is just the
trivial flat metric on $\R^8$, with $r\le0$.  In terms of the
description using $z$ and $v$ introduced in section 3.1, it
corresponds to a degenerate solution at the point $(z,v)=(1,-2)$.
\medskip

\noindent
$\bullet$ 
The solution with $f=3r$ has $Q=0$, and so here we must solve for $b$
using (\ref{absol}).   After making a coordinate transformation
$r\longrightarrow 3r^2/20$, this solution is
%%%%%
\begin{equation}
a^2 = \ft3{10}\, r\,
\Big[1-\Big(\fft{r_0}{r}\Big)^{10/3}\Big]\,,\qquad
b=a\,,\qquad c^2=\ft9{20}\, r^2\,.
\end{equation}
%%%%%
This can be recognised as the previously-known complete metric of
Spin(7) holonomy \cite{brysal,gibpagpop}, as given in
(\ref{spin7metric}).  (The trivial scaling constant $r_0$ arose here
in the integration of the equation for $b$ in (\ref{absol}).)  Note
that because this solution has $f'=$constant, it is not contained within the
general analysis of section 3.1, except as a singular limit that 
corresponds to the point $(z,v)=(-1,+2)$.
\medskip

\noindent
$\bullet$
 The third elementary solution in (\ref{elsol}), $f=r+
r^2/(2\ell^2)$, gives rise to our first examples of 
new complete metrics of Spin(7) holonomy.   After a coordinate
transformation $r\longrightarrow -\ell\, (r+\ell)$, the metric
in local coordinates becomes
%%%%%
\begin{equation}
ds_8^2 = \fft{(r+\ell)^2\, dr^2}{(r+3\ell)(r-\ell)} +
\fft{4\ell^2\, (r+3\ell)(r-\ell)}{(r+\ell)^2}\, R_3^2
+(r+3\ell)(r-\ell)\, (R_1^2+R_2^2) + \ft12(r^2-\ell^2)\,P_a^2\,.
\label{sol1}
\end{equation}
%%%%%
    Assuming that the constant $\ell$ is positive, it is evident that
$r$ should lie in the range $r\ge\ell$.  We can analyse the behaviour
near $r=\ell$ by
defining a new radial coordinate $\rho$, where
$\rho^2=4\ell\,(r-\ell)$.  Near $\rho=0$ the metric approaches
%%%%%
\begin{equation}
ds_8^2 \approx d\rho^2 +\rho^2 \, (R_1^2+ R_2^2 + R_3^2 +
\ft14 P_a^2)\,.
\end{equation}
%%%%%
The quantity $R_1^2+R_2^2+R_3^2 + \ft14 d\Omega_4^2)$
is precisely the metric on the unit
7-sphere, and so we see that near $r=\ell$ the metric $ds_8^2$
smoothly approaches flat $\R^8$.  At large $r$ the function $b$, which
is the radius in the $U(1)$ direction $R_3$, approaches a constant,
and so the metric approaches an $S^1$ bundle over a 7-metric.  This
7-metric is of the form of a cone over $\CP^3$
(described as the $S^2$ bundle over $S^4$) in this asymptotic
region.  The manifold of this new Spin(7) metric, which we are
denoting by $\bA_8$,  is topologically $\R^8$.  In terms of the description
in section 3.1, this solution corresponds to a trajectory in the 
$(z,v)$ plane with $z=1$, and $v$ running from $v=-2$ (at the origin) to 
$v=-\infty$ (in the asymptotic region). Thus we have

\begin{prop} The metric $\bA_8$ given by (\ref{sol1}) with $r>\ell>0$
admits a smooth complete non-singular extension to $\R^8$.
\label{a8prop}
\end{prop}

   We shall use the acronym AC to denote asymptotically conical
manifolds. Thus asymptotically our new metrics behave like a circle
bundle over an AC manifold in which the length of the $U(1)$ fibres
tends to a constant. The acronym ALF is already in use to describe
metrics which tend to a $U(1)$ bundle over an asymptotically Euclidean
or asymptotically locally Euclidean metric with the length of the
fibres tending to a constant. We shall therefore adopt the acronym ALC
to denote manifolds where the base space of the circle bundle is
asymptotically conical.

   Ricci-flat ALC metrics, although not with special holonomy, have
already been encountered.  For example, the higher-dimensional
Taub-NUT metric is defined on $\R^{2n}$ for all $n$ and it is ALC with
the base of the cone being $\CP^{n-1}$. A closely related example is
the Taub-Bolt metric which has the same asymptotics but is defined on
a line bundle over $\CP^{n-1}$. However, the metric on the base of the
cone in this case (with $n=4$) is the Fubini-Study metric on
$\CP^3$, which is quite different from that of the ``squashed''
Einstein metric on $\CP^3$ in our new metrics.  A discussion of ALE
Spin(7) manifolds based on the idea of blowing up orbifolds has been
given in \cite{joyce}.  As far as we are aware, no explicit examples
of this kind have yet been found.

   We get a different complete manifold, which we are denoting by
$\bB_8$, if we take $r$ to be negative.
It is easier to discuss this by instead setting $\ell=-\td\ell$, where
$\td\ell$ and $r$ are taken to be positive.
Thus instead of (\ref{sol1}) we now have
%%%%%
\begin{equation}
ds_8^2 = \fft{(r-\td\ell)^2\, dr^2}{(r-3\td\ell)(r+\td\ell)} +
\fft{4\td\ell^2\, (r-3\td\ell)(r+\td\ell)}{(r-\td\ell)^2}\, R_3^2
+ (r-3\td\ell)(r+\td\ell)\, (R_1^2 +R_2^2) +
\ft12(r^2-\td\ell^2)\, P_a^2\,,
\label{sol2}
\end{equation}
%%%%%
This time, we have $r\ge 3\td\ell$.
Defining $\rho^2= 4\td\ell\, (r-3\td\ell)$,
we find that near $r=3\td\ell$ the metric has the form
%%%%%
\begin{equation}
ds_8^2 \approx d\rho^2 +\rho^2\, (R_1^2+R_2^2+R_3^2) + \td\ell^2\,
P_a^2\,.
\end{equation}
%%%%%

The quantity $(R_1^2+R_2^2+R_3^2)$ is the metric on the unit
3-sphere, and so in this case we find that the metric smoothly
approaches $\R^4\times S^4$ locally, at small distance.  The
large-distance behaviour is the same as for the previous case
(\ref{sol1}).  In the $(z,v)$ plane of section 3.1, this solution 
corresponds to a trajectory with $z=1$, and $v$ running from $v=+2$ (at the 
origin) to $v=+\infty$ (in the asymptotic region).  Thus we have

\begin{prop} The metric $\bB_8$ given by (\ref{sol2}) with $r
>3\td\ell >0$ admits a smooth complete non-singular extension to the
chiral spin bundle over $S^4$.
\label{b8prop}
\end{prop}

   Again we have a complete non-compact ALC metric with Spin(7)
holonomy with the same base.  At short distance, it has the same
structure as the previously-known metric of Spin(7) holonomy, obtained
in \cite{gibpagpop}. 

   We can think of the new manifold $\bA_8$ as providing a smooth
interpolation between Euclidean 8-space at short distance, and an
$S^1$ bundle over ${\cal M}_7$ at large distance, while $\bB_8$
provides an interpolation between the previous Spin(7) manifold of
\cite{brysal,gibpagpop} at short distance and the $S^1$ bundle over
${\cal M}_7$ at large distance.  Here ${\cal M}_7$ denotes the
7-manifold of $G_2$ holonomy that is the $\R^3$ bundle over $S^4$
\cite{brysal,gibpagpop}.

\subsection{General globally-defined solutions}

   Having discussed some special solutions of (\ref{firstorder}) in
section 3.2, and having seen that they include new
complete metrics of Spin(7) holonomy, we now turn to a discussion of
the global structure of the general solutions (\ref{general2}).

    In order to recognise the solutions that give rise to complete 
non-singular metrics, it is helpful first to study the phase-plane
diagram for the first-order equation (\ref{vde}), which can be
expressed as
%%%%%
\begin{equation}
\fft{dz}{d\tau} = 2z\,(1-z^2)\,,\qquad
\fft{dv}{d\tau} = v+2z\,,
\end{equation}
%%%%%
where $\tau$ is an auxiliary ``time'' parameter.   The solutions can
be studied by looking at the flows generated by the 2-vector field
$\{dz/d\tau,dv/d\tau\}=\{ 2z\,(1-z^2),  v+2z\}$ in the $(z,v)$ plane.
For any such flow, it is then necessary to investigate the global
structure of the associated metric (\ref{general2}) for regularity.

   We find that regular solutions can arise in the following four
   cases, namely
%%%%%
\begin{eqnarray}
{\bf (1)}\qquad \bA_8:&& z=1 \ \ \hbox{(fixed)};\qquad v=-2\ \ \hbox{to}\ \
v=-\infty\,,\nn\\
{\bf (2)}\qquad \bB_8 :&& z=1 \ \ \hbox{(fixed)};\qquad v=+2\ \ \hbox{to}\ \
v=+\infty\,,\nn\\
{\bf (3)} \qquad \bB_8^-:&& z_0 \le z\le 1;\qquad v=+2\ \ \hbox{to}\ \
v=+\infty,\qquad (0<z_0 <1)\,,\nn\\
{\bf (4)} \qquad \bB_8^+:&& 1\le z\le z_0;\qquad v=+2\ \ \hbox{to}\ \
v=+\infty,\qquad (1 <z_0 <\infty)\nn\\
&&\qquad\qquad\qquad \quad
\hbox{or}\ z_0 < -1;\ \hbox{(see discussion below)}
\,.\label{sollist}
\end{eqnarray}
%%%%%
Note that $v=\pm\infty$ corresponds to the asymptotic large-distance
region, and in all four cases the metrics have similar asymptotic
structures, precisely as we have already seen in the $\bA_8$ and
$\bB_8$ cases.  The point $v=-2$, $z=1$ corresponds to the
short-distance behaviour of the $\bA_8$ metric, approaching Euclidean
$\R^8$ at the origin where the $S^7$ principal orbits degenerate to a
point.  When $v=2$, on the other hand, we have the short-distance
behaviour seen in the $\bB_8$ metric, approaching $\R^4\times S^4$
locally.  In fact Solution (1) is the metric (\ref{sol1}) on $\bA_8$
found in section 3.2, and Solution (2) is the metric
(\ref{sol2}) on $\bB_8$ found there also.  These both have $k=0$ in
(\ref{vres}).

   Solution (3) arises when $k$ is any positive number, with $z_0$
being the corresponding value of $z$ at which $v(z_0)=2$, with
$0<z_0<1$.  The value of $z_0$ is correlated with the value of $k$,
ranging from $z_0=0$ for $k=\infty$, to $z_0=1$ for
$k=0$.

Near
$z=1$ it follows from (\ref{vres}) that we shall have
%%%%%
\begin{equation}
v=2^{3/4}\, k\, (1-z)^{-1/4}  - 2 + \cdots\,,\qquad
f=c_0\, (1-z)^{-1/2} + \cdots\,,
\end{equation}
%%%%%
where $c_0$ is an arbitrary constant of integration.  Defining
$y\equiv (2c_0)^{-1/2}\, (1-z)^{-1/4}$, we see that as
$z\longrightarrow 1$ we shall have $y\longrightarrow\infty$ and
%%%%%
\begin{equation}
ds_8^2 \approx dy^2 + y^2\, (R_1^2 + R_2^2) + \ft12 y^2\, P_a^2 +
\fft{2\sqrt2\, c_0}{k^2}\, R_3^2\,,
\end{equation}
%%%%%
and so this more general metric has the same large-distance asymptotic
form as do $\bA_8$ and $\bB_8$.  Near $z=z_0$ we shall have $v(z) = 2
+ v'(z_0)\, (z-z_0) +\cdots$, and defining a new radial coordinate $x$
by $(z-z_0)= \ft14 x^2$ near $z=z_0$, we shall have
%%%%%
\begin{equation}
ds_8^2 \approx \fft{f_0}{2z_0\, (1-z_0^2)(1-z_0)\, v'(z_0)}\,
\Big[ dx^2 + v'(z_0)^2\, z_0^2\, (1-z_0)^2\, x^2\,
(R_1^2+R_2^2+R_3^2)\Big]  + f_0\, P_a^2\,,\label{short3}
\end{equation}
%%%%% 
where $f_0$ is the value of $f$ at $z=z_0$.  From (\ref{vde}) we have
that $z_0(1-z_0)\, v'(z_0) =1$, and so we see from (\ref{short3}) that
at short distance the metric (\ref{short3}) approaches $\R^4\times
S^4$ locally.  Thus these more general solution (3) in (\ref{sollist})
with $k>0$ is complete on a manifold that is very similar to the
manifold $\bB_8$ of the solution (\ref{sol2}), with an $S^4$ bolt at
$z=z_0$.  Here, the trajectory in the $(z,v)$ plane runs from
$(z_0,2)$ to $(1,\infty)$.  We shall denote the solution by $\bB_8^-$,
where the superscript indicates that $z$ starts from a value $z_0<1$
at short distance, flowing to $z=1$ asymptotically.  For the case
$k=0$, which leads to the metrics (\ref{sol1}) and (\ref{sol2}), the
quantity $z$ is not a good choice for the radial coordinate, since it
is fixed at $z=1$.  This case can be regarded as a singular limit
within the general formalism we are using here.  Specifically, if we
let $z=1-16 \ep^4\, \td\ell^4\, (r+\td\ell\, )^{-4}$, $k=2^{1/4}\,
\ep$, and choose the integration constant in (\ref{general}) so that
$f=\ft12(r^2-\td\ell^2)$, then upon sending $\ep$ to zero we recover
the metric (\ref{sol2}).

   Solution (4) arises in the case where at large distance $z$ now
approaches 1 from above, and again the flow runs from an $S^4$ bolt at
which $v(z_0)=2$, to the asymptotic region as $z$ approaches 1.  There
are two possibilities, with $z_0$ either being greater than 1, or else
$z_0$ is less than $-1$.  In the latter case $z$ then runs from $z_0$
at the bolt, through $z=-\infty$ to $z=+\infty$, and then down to
$z=1$ in the asymptotic region.  It is useful now to make another
change of radial coordinate, and define $y\equiv 1/z$. This allows the
two regions for $z_0$ to be combined.  The solution for $v$ may now be
written as
%%%%%
\begin{equation}
v= (1-y^2)^{-1/4}\, \Big( \kappa + y\,\, _2F_1[\ft12,\ft34;\ft32; y^2]
\Big)\,.
\end{equation}
%%%%%
The $y$ coordinate then ranges from $y=y_0$ at the bolt to $y=1$ at
infinity, and $-1 \le y_0\le 1$.  The integration constant $\kappa$ is
determined in terms of $y_0$ by the requirement that $v=2$ at $y=y_0$.
It can range between $\kappa=2\sqrt{\pi}\, \Gamma(\ft54)/\Gamma(\ft34)$,
corresponding to $y_0=-1$, and
$\kappa= -2\sqrt{\pi}\, \Gamma(\ft54)/\Gamma(\ft34)$, corresponding
(by taking a suitable limit analogous to the one discussed in footnote
3) to $y_0=+1$.

  A similar analysis to that for Solution 3 above now shows that the
metric in Solution 4 smoothly approaches $\R^4\times S^4$ locally at
$y=y_0$, and that it has the same asymptotic behaviour as the previous
examples.  We shall denote this solution by $\bB_8^+$.  Note that the
simple solution $\bB_8$ in (\ref{sol2}) can be viewed as the
$k\longrightarrow 0$ or $\kappa\longrightarrow -2\sqrt{\pi}\,
\Gamma(\ft54)/\Gamma(\ft34)$ limit of the more complicated $\bB_8^-$
or $\bB_8^+$ solutions respectively.

   We observed at the end of section 2 that a particular example of a
solution of the first-order equations (\ref{firstorder}) is the direct
product metric $ds_8^2=ds_7^2+d\varphi^2$, where $ds_7^2$ is the
Ricci-flat 7-metric of $G_2$ holonomy on the $\R^3$ bundle over $S^4$
\cite{brysal,gibpagpop}, and $\varphi$ is a coordinate on a circle.
We are now in a position to see how this solution can arise as a limit
of our new Spin(7) metrics.  Specifically, it arises as the
$k\longrightarrow \infty$ limit of Solution (3) listed in
(\ref{sollist}).  This is the limit where the constant $z_0$, which
sets the lower limit for the range $z_0\le z\le 1$ for $z$, becomes
zero.  At the same time as sending $k$ to infinity, we can rescale the
fibre coordinate $\varphi$ appearing in $R_3=\ft12 d\varphi+\cdots$,
according to $\varphi\longrightarrow k\, \varphi$.  From (\ref{vres})
and (\ref{general}) it follows that when $k$ becomes very large we
shall have
%%%%%
\begin{equation}
v\longrightarrow \fft{2k\, \sqrt{z}}{(1-z^2)^{1/4}}\,,\qquad
f\longrightarrow \Big(\fft{1+z}{1-z}\Big)^{1/2}\,,\label{decoupling}
\end{equation}
%%%%%
and so in the limit of infinite $k$ the metric (\ref{general2})
becomes
%%%%%
\begin{equation}
ds_8^2 = \fft{dz^2}{4z\, (1-z)^2\, (1-z^2)^{1/2}}
+ \fft{4z}{(1-z^2)^{1/2}}\, (R_1^2+R_2^2)
+\Big(\fft{1+z}{1-z}\Big)^{1/2}\, P_a^2 + d\varphi^2\,.
\end{equation}
%%%%% 
Defining a new radial coordinate $r$ by $r^4=(1+z)\, (1-z)^{-1}$, we
see that this becomes $ds_8^2=ds_7^2 + d\varphi^2$, where
%%%%%
\begin{equation}
ds_7^2 = \fft{2dr^2}{1-r^{-4}} + 2 r^2\, (1-r^{-4})\, (R_1^2 + R_2^2)
+ r^2\, P_a^2\,.
\end{equation}
%%%%%
This can be recognised as the metric of $G_2$ holonomy on the manifold
${\cal M}_7$ of the $\R^3$ bundle over $S^4$, which was constructed in
\cite{brysal,gibpagpop}. (If we had not rescaled the fibre
coordinate $\varphi$ by a factor of $k$ before taking
the limit $k\longrightarrow \infty$, the radius of the $S^1$
would have tended to zero, this limit being taken in the sense of
Gromov-Hausdorff.)  Thus the family of new
Spin(7) manifolds that we are denoting by $\bB_8^-$ has a non-trivial
parameter $k$ such that the $k=\infty$ limit degenerates to ${\cal
M}_7\times S^1$, while the $k=0$ limit reduces to the case $\bB_8$
given by (\ref{sol2}).

\section{Parallel spinor and calibrating 4-form}

    The fact that the metric (\ref{8metans1}), together with
(\ref{firstorder}), has Spin(7) holonomy implies, and is implied by,
the existence of a globally-defined parallel spinor field $\eta$,
satisfying $D\, \eta=0$, where $D\equiv d+ \ft14 \hat\omega_{AB}\,
\Gamma_{AB}$ is the Lorentz-covariant exterior derivative that acts on
spinors in eight dimensions.  Here $\Gamma_{AB} \equiv
\ft12(\Gamma_A\, \Gamma_B-\Gamma_B\, \Gamma_A)$, and $\Gamma_A$ are
the Dirac matrices that generate the Clifford algebra in eight
dimensions.

    After a straightforward calculation, we find that $D$ is given by
%%%%%
\begin{eqnarray}
D&=& d + e^0\, \Big( \fft{\dot c}{c}\Gamma_{08} -\fft{a}{4c^2}\,
\Gamma_{1\hat 1} -\fft{a}{4c^2}\, \Gamma_{2\hat 2} -
\fft{b}{4c^2}\, \Gamma_{3\hat 3}\Big) +
e^1\, \Big( \fft{\dot c}{c}\Gamma_{18} +\fft{a}{4c^2}\,
\Gamma_{0\hat 1} +\fft{a}{4c^2}\, \Gamma_{3\hat 2} -
\fft{b}{4c^2}\, \Gamma_{2\hat 3}\Big) \nn\\
&&+ e^2\, \Big( \fft{\dot c}{c}\Gamma_{28} +\fft{a}{4c^2}\,
\Gamma_{0\hat 2} -\fft{a}{4c^2}\, \Gamma_{3\hat 1} +
\fft{b}{4c^2}\, \Gamma_{1\hat 3}\Big) +
e^3\, \Big( \fft{\dot c}{c}\Gamma_{38} +\fft{a}{4c^2}\,
\Gamma_{2\hat 1} -\fft{a}{4c^2}\, \Gamma_{1\hat 2} -
\fft{b}{4c^2}\, \Gamma_{0\hat 3}\Big) \nn\\
&&+ e^{\hat 1}\, \Big(\fft{\dot a}{2a}\, \Gamma_{\hat 1 8}+
 (\fft{a}{4c^2} -\fft1{4a})\,
(\Gamma_{01}+\Gamma_{23})-\fft{b}{4a^2}\, \Gamma_{\hat 2\hat 3}\Big)
\nn\\
&&+ e^{\hat 2}\, \Big(\fft{\dot a}{2a}\, \Gamma_{\hat 2 8}+
 (\fft{a}{4c^2} -\fft1{4a})\,
(\Gamma_{02}+\Gamma_{31})-\fft{b}{4a^2}\, \Gamma_{\hat 3\hat 1}\Big)
\nn\\
&&+ e^{\hat 3}\, \Big(+\fft{\dot b}{2b}\, \Gamma_{\hat 3 8}+
 (\fft{b}{4c^2} -\fft1{4b})\,
(\Gamma_{03}+\Gamma_{12})+(\fft{b}{4a^2}-\fft1{2b})\, \Gamma_{\hat1\hat2}
\, \Big)\,.
\end{eqnarray}
%%%%%
Substituting in the first-order equations (\ref{firstorder}), 
this becomes
%%%%%
\begin{eqnarray}
D &=& d + e^0\, \Big(\fft{a}{4 c^2}\, (2\Gamma_{08} -\Gamma_{1\hat 1} -
\Gamma_{2\hat 2}) + \fft{b}{4c^2}\, (\Gamma_{08} -\Gamma_{3\hat
3})\Big) \nn\\
&& + e^1\, \Big(\fft{a}{4 c^2}\, (2\Gamma_{18} +\Gamma_{0\hat 1} +
\Gamma_{3\hat 2}) + \fft{b}{4c^2}\, (\Gamma_{18} -\Gamma_{2\hat
3})\Big)\nn\\
&&
+ e^2\, \Big(\fft{a}{4 c^2}\, (2\Gamma_{28} +\Gamma_{0\hat 2} -
\Gamma_{3\hat 1}) + \fft{b}{4c^2}\, (\Gamma_{28} +\Gamma_{1\hat
3})\Big)\nn\\
&&+ e^3\,\Big(\fft{a}{4 c^2}\, (2\Gamma_{38} -\Gamma_{1\hat 2} +
\Gamma_{2\hat 1}) + \fft{b}{4c^2}\, (\Gamma_{38} +\Gamma_{0\hat 3})
\Big)\nn\\
&& + e^{\hat 1}\, \Big( (\fft1{4a}-\fft{a}{4c^2})\, (2\Gamma_{\hat1 8} -
\Gamma_{01} -\Gamma_{23}) -\fft{b}{4a^2}\, (\Gamma_{\hat1 8} +
\Gamma_{\hat2\hat3} )\Big) \nn\\
&&+ e^{\hat 2}\, \Big( (\fft1{4a}-\fft{a}{4c^2})\, (2\Gamma_{\hat2 8} -
\Gamma_{02} -\Gamma_{31}) -\fft{b}{4a^2}\, (\Gamma_{\hat2 8} +
\Gamma_{\hat3\hat1} )\Big) \nn\\
&&+ e^{\hat 3}\, \Big(\fft{b}{4a^2}\, (\Gamma_{\hat3 8} +
\Gamma_{\hat1 \hat2}) - \fft{b}{4c^2}\, (2\Gamma_{\hat3 8}
-\Gamma_{03} - \Gamma_{12}) -\fft1{4b}\, (2\Gamma_{\hat1 \hat2} +
\Gamma_{03} + \Gamma_{12})\Big)\,.
\end{eqnarray}
%%%%%
It is now straightforward to see that in this frame, a spinor $\eta$
satisfies $D\, \eta=0$ if it has constant components, and if in
addition it satisfies the projection conditions
%%%%%
\begin{equation}
(2\Gamma_{08} - \Gamma_{1\hat 1}-\Gamma_{2\hat 3})\, \eta=0\,,
\quad
(\Gamma_{08} -\Gamma_{3\hat 3})\, \eta=0\,,\quad
(2\Gamma_{18} + \Gamma_{0\hat 1} + \Gamma_{3\hat 2})\,\eta=0\,.
\label{gammaproj}
\end{equation}
%%%%% 
These conditions define a unique spinor, up to overall scale, thus providing
another proof that the metrics have Spin(7) holonomy.

    The covariantly-constant self-dual 4-form $\Phi$
given in (\ref{calib}), known as the
Cayley form, provides a calibration of the Spin(7) manifold.  Thus we have
%%%%%
\begin{equation}
|\Phi(X_1,X_2,X_3,X_4)| \le 1\,,
\end{equation}
%%%%%
where $(X_1,X_2,X_3,X_4)$ denotes any quadruple of orthonormal
vectors.  This can be seen from (\ref{calib}), or else by noting
that the components of $\Phi$ can be expressed in terms of the
parallel spinor $\eta$ as
$\Phi_{ABCD} = \bar\eta\, \Gamma_{ABCD}\, \eta$.  A
calibrated submanifold, or Cayley submanifold, $\Sigma$, is one where
for each point of $\Sigma$
%%%%%
\begin{equation}
|\Phi(X_1,X_2,X_3,X_4)|=1\,,
\end{equation}
%%%%%
where the orthonormal vectors $X_i$ are everywhere tangent to
$\Sigma$.  By inspecting (\ref{calib}) we therefore see that the
$S^4$ zero section of the bundle of chiral spinors is a Cayley
submanifold, and hence it is volume minimising in its homology class.

    Thus in summary we have

\begin{prop} The gradient flow equations (\ref{firstorder}) are the
necessary and sufficient conditions for the local metric (\ref{8metans1}) to
have holonomy Spin(7).  The covariantly-constant spinor is defined by 
(\ref{gammaproj}), and the Cayley 4-form is given by (\ref{calib}).
\label{spin7prop}
\end{prop}

\section{$L^2$-normalisable harmonic 4-forms in $\bA_8$ and $\bB_8$}

     In this section, derive the equations for harmonic 4-forms in the 
Spin(7) manifolds.  We obtain explicit $L^2$ normalisable harmonic 4-forms
for each of the new Spin(7) 8-manifolds $\bA_8$ and $\bB_8$.
Specifically, we obtain one such 4-form, which is anti-self-dual, for the
manifold $\bA_8$ that is topologically $\R^8$, and two such 4-forms, one
of each duality, for the manifold $\bB_8$ of the chiral spin bundle over
$S^4$.

     The structure of the harmonic 4-forms turns out to be closely
related to that of the calibrating 4-form $\Phi$ given in (\ref{calib}).
Thus we define
%%%%%
\begin{equation}
G_\4^\pm = \omega_\4\pm {*\omega_\4}\,,\label{g4def}
\end{equation}
%%%%%
where
%%%%%
\begin{equation}
\omega_\4 \equiv u_1\,e^0\wedge e^1\wedge e^2\wedge e^3
- u_2\, (e^{\hat 2}\wedge e^{\hat 3}\wedge
\hat J^1 + e^{\hat 3}\wedge e^{\hat 1}\wedge \hat J^2) +
u_3\, e^{\hat 1}\wedge e^{\hat 2}\wedge \hat J^3\,.\label{omegadef}
\end{equation}
%%%%% 
$G^\pm_\4$ will be harmonic if $dG^\pm_\4=0$.  This implies that
%%%%%
\begin{eqnarray}
\pm\fft{d(c^4\, u_1)}{dt} - 2b\, c^2\, u_2 + 4  a\, c^2\,
u_3 &=&0\,,\nn\\
\pm\fft{d(a^2\, c^2\, u_2)}{dt} - a^2\, b\, u_1 +
b\, c^2\,  u_2  + 2 a\, c^2\, u_3 &=&0\,,\label{4formeq}\\
\pm\fft{d(a\, b\, c^2\, u_3)}{dt} + a^2\, b\,  u_1 + b\, c^2\, u_2
&=&0\,.\nn
\end{eqnarray}
%%%%%%
where the $\pm$ signs correspond to self-dual and anti-self-dual
respectively.  (The Cayley form given in (\ref{calib}) is a particular
solution, corresponding to taking $u_1=u_2=-1$, $u_3=1$.)

    In the case of the new Spin(7) manifolds $\bA_8$ and $\bB_8$, with
their simple metrics (\ref{sol1}) and (\ref{sol2}), we can now obtain
explicit results for $L^2$ harmonic 4-forms.  In the remainder of this
section, we shall for convenience set the scale parameters $\ell$ and
$\td\ell$ in the metrics (\ref{sol1}) and (\ref{sol2}) to
unity.   Care must be exercised when taking the square roots of
$a^2$, $b^2$ and $c^2$ in the metrics (\ref{sol1}) and (\ref{sol2}),
if one wants the functions $a$, $b$ and $c$ to solve precisely the
first-order equations (\ref{firstorder}), since these equations are
sensitive to the signs of $a$, $b$ and $c$.  (Of course there are
equivalent first-order equations that differ by precisely these sign
factors, and which also imply solutions of the Einstein equations.)
We are assuming here that the signs are chosen so that precisely
(\ref{firstorder}) are satisfied.  This can be achieved by taking all
square roots to be positive, except for $b$ in the case of
(\ref{sol1}) on $\bA_8$.

     For the metric (\ref{sol1}) on the manifold $\bA_8$ that is
topologically $\R^8$, we find
that there is a normalisable harmonic 4-form that is anti-self-dual,
{\it i.e.}, the lower choice of the sign is used in (\ref{g4def})
and (\ref{4formeq}).  The solution is given by
%%%%%%%
\begin{equation}
u_1 =\fft{2}{(r+1)^3(r+3)}\,,\quad
u_2 =-\fft{r^2 +10r + 13}{(r+1)^3(r+3)^3}\,,\quad
u_3 =-\fft2{(r+1)^2(r+3)^3}\,.\label{aasd}
\end{equation}
%%%%%%
The norm of the harmonic anti-self-dual 4-form is then given by
%%%%%
\begin{equation} 
|G_\4|^2 = 48(u_1^2 +2 u_2^2 + 4 u_3^2)=\fft{96(3r^4+44r^3
+ 242r^2 + 492r +339)}{(r+1)^6(r+3)^6}\,.\label{g4square1anti}
\end{equation} 
%%%%%
Clearly $G_\4$ is $L^2$-normalisable, and in fact we have
$\int_1^\infty \sqrt{g}\, |G_\4|^2\, dr = 9/4$.  We have chosen
the integration constants from (\ref{4formeq}) appropriately in order
to select the solution in $L^2$.  (There also exists a solution for 
a self-dual harmonic 4-form.  It can be made square integrable at small
distance, but there is no choice of integration constants for which it
is $L^2$ normalisable, owing to its large distance behaviour.)

      For the metric (\ref{sol2}) on $\bB_8$, the bundle of chiral spinors
over $S^4$, we find that there exists a normalisable harmonic 4-form
that is anti-self-dual, {\it i.e.}, the lower choice of sign is used
in (\ref{g4def}) and (\ref{4formeq}).  The solution is given by
%%%%%
\begin{eqnarray}
&&u_1=\fft{2(r^4+8r^3+34r^2-48r+21)}{(r-1)^3(r+1)^5}\,,\qquad
u_2= - \fft{r^4+4r^3-18r^2+52r-23}{(r-1)^3(r+1)^5}\,,\nn\\
&&u_3=\fft{2(r^2+14r-11)}{(r-1)^2(r+1)^5}\,.\label{basd}
\end{eqnarray}
%%%%%
The square of the anti-self-dual 4-form is given by
%%%%%
\begin{equation}
|G_\4|^2 = \fft{96(3r^8 +40r^7 + 252r^6 +1064r^5+2506r^4
-12936r^3 +18284r^2 -10824r +2379)}{(r-1)^6(r+1)^{10}}\,,
\label{g4square2anti}
\end{equation}
%%%%%
and its $L^2$-normalisability can be seen by noting that
$\int_3^\infty  \sqrt{g}\, |G_\4|^2\, dr = 189/16$.

   We also find a second $L^2$-normalisable harmonic 4-form in the
new Spin(7) manifold $\bB_8$. This 4-form is self-dual, with the upper
sign chosen in (\ref{g4def}) and (\ref{4formeq}) and is given by 
%%%%%
\begin{eqnarray}
&&u_1= -\fft{2(5r^3-9r^2+15r-3)}{(r-1)^3\, (r+1)^4}\,,\nn\\
&&u_2 = \fft{(r-3)(5r^2-2r+1)}{(r-1)^3\, (r+1)^4}\,,\quad
u_3=-\fft{2(r-3)}{(r-1)^2\, (r+1)^4}\,.\label{bsd}
\end{eqnarray}
%%%%%
In contrast to the previous harmonic 4-forms, there is no linear relation 
between the functions $u_1$, $u_2$ and $u_3$ here.  The magnitude of $G_\4$
is given by
%%%%%
\begin{equation}
|G_\4|^2 = \fft{96(75r^6-350r^5+829 r^4-932 r^3+885 r^2-414 r+99)}{
               (r-1)^6\, (r+1)^8}\,.\label{g4square2self}
\end{equation}
%%%%%
It integrates to give $\int_3^\infty \sqrt{g}\, |G_\4|^2\, dr =189/4$.

   It is interesting to note that for the anti-self-dual harmonic
4-form on $\bA_8$, given by (\ref{aasd}), we can write it in terms of
a globally-defined potential, $G_\4=dB_\3$.  Specifically, we find
that $B_\3$ can be written as 
%%%%%
\begin{equation}
B_\3= -(r-1)^2\, \Big[ \fft1{(r+1)^2}\, R_1\wedge R_2\wedge R_3 
+\fft{1}{8(r+3)^2}\,(R_1\wedge J^1 + R_2\wedge J^2)+
\fft{(r+5)}{4(r+1)(r+3)^2}\, R_3\wedge J^3 
\Big]\,.\label{b3exp}
\end{equation}
%%%%%
One can see from (\ref{sol1}) that this has a vanishing magnitude
$|B_\3|^2$ at $r=1$.  On the other hand the analogous expressions for
the potential $B_\3$ for the two harmonic 4-forms (\ref{basd}) and
(\ref{bsd}), which are similarly expressible as functions of $r$ times
the three 3-form structures in (\ref{b3exp}), turn out to have a
diverging magnitude at $r=3$.  In all three cases the $r$-dependent
prefactors tend to constants at infinity.

   Our results on harmonic forms are summarised in the following

\begin{prop} The metric $\bA_8$ of proposition (\ref{a8prop}) has
Spin(7) holonomy and admits an $L^2$ harmonic 4-form, given by
(\ref{g4def}) and (\ref{aasd}), whose duality is
opposite to that of of the Cayley form.  The metric $\bB_8$ of
proposition (\ref{b8prop}) has Spin(7) holonomy and admits both an
anti-self-dual $L^2$ harmonic 4-form, given by (\ref{g4def}) and
(\ref{basd}), and a self-dual $L^2$ harmonic 4-form, given by
(\ref{g4def}) and (\ref{bsd}).
\label{harmprop}
\end{prop}

\section*{Acknowledgements}

   We should like to thank Michael Atiyah, Andrew Dancer, Nigel
Hitchin and Dominic Joyce for useful discussions.
H.L.~and C.N.P.~are grateful to UPenn for hospitality and financial
support during the course of this work.  M.C.~is supported in part by
DOE grant DE-FG02-95ER40893 and NATO grant 976951; H.L.~is supported
in full by DOE grant DE-FG02-95ER40899; C.N.P.~is supported in part by
DOE DE-FG03-95ER40917.  G.W.G.~acknowledges partial support from PPARC
through SPG\#613.

\end{document}